# Finitary Codes, a short survey

## Jacek Serafin[1]

*Wrocław University of Technology*

**Abstract:** In this note we recall the importance of the notion of a finitary isomorphism in the classification problem of dynamical systems.

## 1. Introduction

The notion of isomorphism of dynamical systems has been an object of an extensive research from the very beginning of the modern study of dynamics. Various invariants of isomorphism have been introduced, like ergodicity, mixing of different kinds, or entropy. Still, the task of deciding whether two dynamical systems are isomorphic remains a notoriously difficult one.

Independent processes, referred here to as Bernoulli schemes (defined below), are most easily described examples of dynamical systems. However, even in this simplest situation, the isomorphism problem remained a mystery from around 1930, when it was stated, until late 1950's, when Kolmogorov implemented Shannon's ideas and introduced mean entropy into ergodic theory; later Sinai was able to compute entropy for Bernoulli schemes and show that Bernoulli schemes of different entropy were not isomorphic.

In his ground-breaking work Ornstein (see [20]) produced metric isomorphism between Bernoulli schemes of the same entropy, and later he and his co-workers discovered a number of conditions which turned out to be very useful in studying the isomorphism problem. Among those conditions are notions of Weak Bernoulli, Very Weak Bernoulli and Finitely Determined. As those will not be of our interest in this note, we shall not go into the details, and the interested reader should for example consult [20] or [24]. What is relevant for us here is the fact that the mappings which Ornstein constructed had the following disadvantageous property: in order to determine a target sequence one should examine the entire input sequence (the whole infinite past and future). On the other hand, it seemed necessary for applications that the codes should enjoy some kind of continuity property.

## 2. Finitary Codes

Numerous attempts had been done, after Kolmogorov introduced entropy into the study of dynamical systems, to construct effective codes between independent processes and between Markov processes. Meshalkin ([17]), Blum and Hanson [3], Monroy and Russo [18], among others, have successfully constructed isomorphisms between specific independent or Markov processes, and all of those codes enjoyed the following property: to determine a value of any given coordinate in an output sequence, one should only examine finitely many coordinates of the source sequence,

---

this finite number of coordinates depending upon the input sequence under consideration.

A fundamental object of our study is a *dynamical system*: a Lebesgue probability space $(X, \mathcal{A}, \mu)$ together with a measure-preserving, ergodic automorphism $S$ of $X$. However, as it is well-known ([14]) that every ergodic dynamical system of finite entropy can be represented as a shift transformation on a sequence space with a finite underlying alphabet, from now on we shall focus our attention to the case of sequence spaces. Consequently, we assume that $X$ is a sequence space, $\mathcal{A}$ is a product $\sigma$-algebra generated by the coordinate mappings, $S$ is the left shift transformation and $\mu$ is a shift-invariant measure. We are now ready to define the principal objects of our interest.

**Definition 1.** A *homomorphism* (factor map) $\phi$ from a dynamical system $(X, \mathcal{A}, \mu, \mathcal{S})$ to a dynamical system $(Y, \mathcal{B}, \nu, \mathcal{T})$ is a measurable map $\phi$ from a subset of measure one of $X$ to $Y$ such that $\nu = \mu \phi^{-1}$ and $\phi S = T \phi$. The homomorphism $\phi$ is called *finitary* if it becomes a continuous map after discarding subsets of measure zero from both spaces $X$ and $Y$.

If $\phi$ is invertible and $\phi^{-1}$ is continuous outside subsets of measure zero, then we call $\phi$ a finitary isomorphism.

For shift spaces, an equivalent definition is the following:

**Definition 2.** A homomorphism $\phi$ between $X$ and $Y$ is finitary if for almost every $x \in X$ there exist positive integers $q = q(x), r = r(x), q \leq r$, such that if $y \in X$ and $[x_{-q}, \ldots, x_r] = [y_{-q}, \ldots, y_r]$ and if $\phi(y)$ is defined, then $(\phi(x))_0 = (\phi(y))_0$.

Let us note in passing that definitions of almost-continuity can be given in a more general setting, and for this we refer the interested reader to [5].

The following random variable, called *code length* function, will be of special interest to us:

**Definition 3 (code length).** For $x \in X$, let $q = q(x)$ and $r = r(x)$ be minimal positive integers, as in the previous definition. The code length for $x \in X$ is $C(x) = q(x) + r(x)$.

Random variables $q$ and $r$ are sometimes called *memory* and *anticipation* of the coding, respectively.

**Definition 4.** We say that dynamical systems are finitarily isomorphic with finite expected code times (fect) if both the finitary isomorphism $\phi$ and its inverse $\phi^{-1}$ have finite expected code lengths.

A systematic study of finitary coding had begun with the works of Keane and Smorodinsky ([10], [11], [12]), and Denker and Keane ([5], [6]). In their 1977 paper [10], and in subsequent articles [11], [12], Keane and Smorodinsky developed a theory and methodology of finitary coding, creating a new area of research, which has been (and is still being) extended in a multitude of ways. In this note we want to recall the marker methods of [10], and later discuss the existence of finitary homomorphism with finite expected code length or a finitary isomorphism with fect.

It is perhaps worth noticing here that, despite major developments in the field, some of the most fundamental questions regarding classification remain open for more than 20 years now.

## 2.1. Case of different entropies

In this section we recall the main ideas behind the Keane-Smorodinsky construction [10]. The basic object of our study is a space of doubly-infinite sequences drawn from a finite alphabet consisting of $a$ symbols ($a \geq 2$), $X = \{1, \ldots, a\}^{\mathbb{Z}}$, equipped with the product $\sigma$-algebra $\mathcal{A}$, product measure $\mu = \mathbf{p}^{\mathbb{Z}}$ and the left shift transformation $S$. Here $\mathbf{p} = (p_1, \ldots, p_a)$ is a strictly positive probability vector assigning probabilities to symbols $1, \ldots, a$. A quadruple

$$(X, \mathcal{A}, \mu, \mathcal{S})$$

is commonly referred to as a *Bernoulli scheme* based on a probability vector $\mathbf{p}$, and will be denoted $BS(p)$. It is well-known that the entropy of $BS(p)$ equals $h = h(p) = -\sum_i p_i \log p_i$. We will need another Bernoulli scheme, $BS(q)$, based on a probability vector $\mathbf{q} = (q_1, \ldots, q_b)$ on $b$ symbols, here the shift space is $\bar{X} = \{1, \ldots, b\}^{\mathbb{Z}}$, the product $\sigma$-algebra is $\mathcal{B}$, product measure is $\nu = \mathbf{q}^{\mathbb{Z}}$ and the left shift transformation is $T$, entropy of $T$ is $\bar{h} = h(q)$. In 1969 Ornstein (see [20]) proved that Bernoulli shifts of the same entropies were isomorphic and also showed the following:

**Theorem 1.** *If $h > \bar{h}$ then Bernoulli scheme $BS(q)$ is a homomorphic image of $BS(p)$.*

Later ([10]) Keane and Smorodinsky strengthened both the unequal and equal entropies statements to the case of finitary coding. In this section we want to focus on the following statement:

**Theorem 2 ([10]).** *If $h > \bar{h}$, then there exists a finitary homomorphism from $BS(p)$ to $BS(q)$.*

Before we continue, let us mention that the isomorphism result [11] relies upon a beautiful refinement and improvement of the methods developed for unequal entropies case. An excellent exposition of the finitary isomorphism result (which has become standard in ergodic theory) appeared in Petersen's book [24] (see also [4] or [23]).

We now proceed to describe the Keane–Smorodinsky construction (see also [1]) in the case of unequal entropies, and then recall some subsequent results in which the new techniques were applied.

Basic reduction allows to assume that there are two blocks, called *markers*, one in each scheme, of the same length and the same probability of appearance; it is also quite natural to demand that the coding procedure should map markers to markers, the main difficulty is in inventing a code for blocks, called *fillers*, occurring between markers. We define a marker, for either scheme, to be a block

$$M = 1^{k-1}2 = \underbrace{1 \ldots 1}_{(k-1)-times} 2,$$

consisting of $k-1$ consecutive 1's followed by a 2. Let us note that a marker has the following *non-overlapping* property: none of its initial subblocks is equal to its terminal subblocks of the same length; this property guarantees the shift-invariance of the coding.

By ergodicity, almost every source sequence $x$ in $X$ splits into *runs* of markers labeled in a natural manner by $\pm 1, \pm 2, \ldots$ and *separating blocks* labeled $0, \pm 1, \pm 2, \ldots$ We assume that the 0-th coordinate of $x$ is covered by either the run of markers labeled −1 or the subsequent 0-th separating block. By $u_j$ we denote the number of markers in the $j$-th run while $l_j$ stands for the length of the $j$-th separating block.

For every $r = 1, 2, \ldots$ we denote by $s_r = s_r(x)$ the *skeleton of rank $r$*. This is defined as the truncation of $x$ to a finite segment around 0 such that the separating blocks in $x$ are replaced by gaps of the same length, and with the property that the extreme left and right runs of markers contain each at least $r$ markers while the internal runs, if any, contain each less than $r$ markers. Moreover, neither the immediately preceding nor the immediately following $k$-block of $x$ is a marker block. We denote by $-m_r < 0$ and $n_r > 0$ the label of the first and the last run of markers in $s_r$, respectively. We may draw the following picture of the skeleton $s_r(x)$ of rank $r$ at $x$, as a sequence of markers and spaces between the markers:

$$M^{u_{-m}} \underset{l_{-(m-1)}}{\rule{1cm}{0.4pt}} M^{u_{-(m-1)}} \cdots \underset{l_{-1}}{\rule{1cm}{0.4pt}} M^{u_{-1}} \underset{l_0}{\rule{1cm}{0.4pt}} M^{u_1} \underset{l_1}{\rule{1cm}{0.4pt}} \cdots M^{u_{n-1}} \underset{l_{n-1}}{\rule{1cm}{0.4pt}} M^{u_n}$$

where $m = m(r), n = n(r) \geq 1$, $l_i \geq 1$, for $i = -(m-1), \ldots, n-1$ and $u_{-m}, u_n \geq r > u_{-(m-1)}, \ldots, u_{n-1}$.

Clearly the rank one skeletons consist of two marker runs separated by one filler block. In order to avoid ambiguity we assume that the zero coordinate of $x$ appears in the 'interior' of $s_1(x)$, i.e. it corresponds to the blank part of $s_1$. The skeleton of rank $r$ is obtained by looking to the left and to the right for the first appearance of $M^r$. For a skeleton $s$ we denote by $l(s)$ the length of $s$ minus the last run of markers,

$$l(s_r) = ku_{-m_r} + l_{-m_r+1} + \cdots + ku_{-1} + l_0 + ku_1 + \cdots + l_{n_r-1},$$

as the final block of markers is only needed to determine the occurrence of $s_r(x)$ but is not considered to be part of that occurrence. A block in $x$ occurring along a single run of markers followed by a separating block will be called an *order one filler*. The concatenation of all the order one fillers in $s_r$ will be referred to as the *filler of $s_r$*. Clearly the length of the filler is equal to $l(s_r)$.

For a fixed non-indexed skeleton $s$ the *filler measure $\mu_s$* is defined on the $l(s)$-blocks as the projection of the conditional measure $\mu(\cdot|\mathcal{S})$ where $\mathcal{S}$ is the event that $s$ occurs at $[0, l(s) - 1]$ in $x$. According to [10], the filler measure $\mu_s$ is the product of the filler measures corresponding to order one subskeletons of $s$. Regardless of the skeleton rank, filler measures will be denoted by $\mu_s$.

Recall that the *marker process* is a stationary 0-1 process $\hat{X}$ defined by $\hat{x}_i = 0$ iff $x_i \ldots x_{i+k-1}$ is a marker block. If the marker length $k$ is sufficiently large then the entropy of the marker process can be made as small as needed. It is clear that both Bernoulli schemes have the marker process as a common factor, and that the skeleton structure at $x$ only depends upon the marker process $\hat{x}$. We define the *filler entropy* $f = h(X) - h(\hat{X})$.

Let us fix $\epsilon < (f - \bar{h})/3$. A filler $F$ in the skeleton $s_r(x)$ of the source sequence $x \in X$, is called *good* if $\mu_s(F) \leq e^{-l(s_r)(f-\epsilon)}$. On the other hand, a corresponding $l(s_r)$-block $\bar{F}$ in $\bar{x} \in \bar{X}$ will be called a good filler if $\bar{\mu}(\bar{F}) \geq e^{-l(s_r)(\bar{h}+\epsilon)}$. According to [10], only good fillers will be encoded to good fillers. If a filler is bad, it will be encoded as a part of a longer good filler at a later stage. The coding is carried out for a given source sequence $x$ by looking at the ascending skeletons $s_r(x)$, $r = 1, 2, \ldots$ By means of an "assignment" defined in [10], the filled skeleton $s_r(x)$ will be encoded in a consistent way if the filler $F$ is good, except for a small set of exceptional cases. Let us now be more specific about the coding procedure.

The set of all possible fillers $\bar{\mathcal{F}}(s_r)$ for the low entropy scheme is divided into equivalence classes. If $\bar{F} \in \bar{\mathcal{F}}(s_r)$ is good then no other filler in $\bar{\mathcal{F}}(s_r)$ is equivalent to it. Let us suppose that $\bar{F}$ is bad; it is then possible that the restriction of $\bar{F}$ to a subskeleton $s'$ of $s_r$ is a good filler for that subskeleton $s'$. The equivalence class of a bad filler $\bar{F}$ consists of all bad fillers $\bar{G}$ with the property that $\bar{G}$ has the same collection (as $\bar{F}$) of subskeletons on which its restrictions are good, and those good restrictions agree with the corresponding (good) restrictions of $\bar{F}$. Let us write $\bar{F} \sim \bar{G}$ if $\bar{F}$ and $\bar{G}$ are equivalent. A simple combinatorial argument in [10] guarantees that there are no more than $2^{m_r+n_r-1+(\bar{h}+\epsilon)l(s_r)}$ equivalence classes in $\bar{\mathcal{F}}(s_r)$; we shall need that estimate later. An important step toward defining the desired code is a notion of *partial assignments*. For a fixed skeleton $s$, a partial assignment $P_s$ assigns to each element $\bar{F}$ of $\bar{\mathcal{F}}(s)$ a subset $P_s(\bar{F})$ of $\mathcal{F}(s)$, in such a way that

$$\nu_s(\bar{F}) \leq \mu_s(P_s(\bar{F})).$$

Partial assignment is an example of a more general notion of a *society* (for a detailed discussion of societies see [10] or [24]). The partial assignment $P_s$ is *good*, if it respects the equivalence classes of $\bar{F}(s)$, i.e. if

$$\bar{F} \sim \bar{G} \Rightarrow P_s(\bar{F}) = P_s(\bar{G}).$$

If we suppose that $\bar{F} \sim \bar{G}$, then the above condition guarantees that each filler $F$ in $P_s(\bar{F})$ will be assigned to $\bar{F}$ and $\bar{G}$ and all other elements of the equivalence class of $\bar{F}$. The Shannon–McMillan–Breiman theorem implies, however, that at some finite stage $\bar{F}$ becomes a part of a longer good filler $\bar{H} \in \bar{\mathcal{F}}(s_r)$ for a skeleton $s_r$ which has $s$ as a subskeleton. Since $\bar{H}$ is good, no other element of $\bar{\mathcal{F}}(s_r)$ is equivalent to it, and each $F \in P_{s_r}(\bar{H})$ is assigned to exactly one filler in $\bar{\mathcal{F}}(s_r)$, that filler being of course $\bar{H}$. Keane and Smorodinsky show in [10], using a version of the marriage lemma, that partial assignments can be consistently extended to so-called global assignments, in such a way that if at a finite stage of the above procedure a filler $F \in \mathcal{F}(s)$ is assigned to the unique $\bar{F} \in \bar{\mathcal{F}}(s)$, then the assignments which take place for skeletons for which $s$ is a subskeleton, respect the $F \mapsto \bar{F}$ assignment. Finally, $\bar{F}$ is defined as a homomorphic image of $F$. A natural question which one might ask, having constructed a finitary coding, is whether the average code length is finite? Even though it was believed that the expected code length should be finite in the case $h(p) > h(q)$, the proof of this statement, quite nontrivial, was given much later in [30]. Shortly after [10], Akcoglu, del Junco and Rahe extended the result of Keane and Smorodinsky. They constructed a finitary coding between an ergodic Markov shift $X$ and a mixing Markov shift $Y$ of smaller entropy. Their construction is quite similar to that of [10], the essential role is being played by the low entropy marker process. Informally speaking, the presence of markers makes it possible to represent almost every source sequence $x \in X$ as an ascending nested family of words, which fill longer and longer skeletons determined by the marker process. Fillers of sufficiently large rank are encoded to corresponding fillers in $Y$ thus eventually defining the required finitary coding $\phi : X \to Y$. Akcoglu, del Junco, and Rahe also claim without proof that the code length should have a finite expectation. We will now very shortly sketch the proof of the main result in [8], which is a significant improvement of the result in [30], that the expected code length between Markov processes of unequal entropies is finite; we also indicate a number of differences between [1] and [10], as we proceed.

One of the differences between [10] and [1] lies in the choice of markers. Unlike in [10], where the authors used $1^{k-1}2$ as a marker, in the Markov case [1] a *marker* $M$ is a collection of blocks of the same length $k$ such that each word in $M$ begins with the same symbol $a_1$, no word in $M$ overlaps a word in $M$ and arbitrarily long concatenations of words from $M$ occur with positive probabilities. Also, the length $k$ can be chosen arbitrarily large and the probability that a marker occurs at a given position decays exponentially with $k$.

The coding between two mixing Markov processes is achieved in two steps, as a composition of two codes, using an intermediate Bernoulli scheme. In the first step, referred to as Markov-to-Bernoulli coding, we study a mixing Markov process $(X, \mu, T)$ and a Bernoulli process $(\bar{X}, \bar{\mu}, \bar{T})$ with $h(\bar{X}) = \bar{h} < h = h(X)$. A marker in $X$ can be selected as a single word $a_1 \ldots a_k$ in such a way that the filler entropy $f$ still exceeds the entropy $\bar{h}$; no marker is needed in the Bernoulli process $\bar{X}$. Good and bad fillers are defined similarly as in the Bernoulli case, and the coding procedure follows to a large extent that of [10].

As far as the code length is concerned, the main results of [8] are the following:

**Theorem 3.** *Let the processes $X$, $\bar{X}$ be mixing Markov and Bernoulli (or Bernoulli and Markov), respectively. If $h(X) > h(\bar{X})$ then there exists a finitary coding from $X$ to $\bar{X}$ such that for every $p < 2$ the code length is an $L^p$ random variable.*

**Lemma 1.** *Let $X, Y, Z$ be arbitrary stationary processes and let $\phi : X \to Y$ and $\psi : Y \to Z$ be finitary codes. Assume that for some $p_1, p_2 > 1$ with $p_1 \leq p_2 + 1$ the code length of $\phi$ is in $L^p$ for all $p < p_1$ and the code length of $\psi$ is in $L^p$ for all $p < p_2$. Then the composed code $\psi \circ \phi$ has code length in $L^p$ for all $p < p_1 p_2/(p_2+1)$.*

Composing the Markov-to-Bernoulli with Bernoulli-to-Markov codes, we obtain:

**Theorem 4.** *Let $X_1$ and $X_2$ be mixing Markov processes such that $h(X_1) > h(X_2)$. Then there exists a finitary coding from $X_1$ to $X_2$ such that the code length is in $L^p$ for all $p < 4/3$.*

Let us recall a basic fact about the nature of the Keane-Smorodinsky coding, that markers are mapped onto markers and fillers onto fillers, so the contents of good skeletons in one scheme define the contents of appropriate skeletons in the other scheme. Consequently, the code length function $C$ only takes values $c_{r+1} = k + l(s_{r+1}) + k(r + 1)$ (depending upon the marker process), which are skeleton lengths plus the length of $r + 1$ markers in the terminal marker occurrence plus $k$ which stands for a number of entries which have to be examined to make sure that there is no marker preceding the initial run of markers. From a combinatorial bound on the number of equivalence classes in $\bar{\mathcal{F}}(s_r)$ it follows that the conditional probability, given a marker structure of $x$, that $C(x) \geq c_{r+1}$, is bounded by (see [10], lemma 14)

$$2^{m_r+n_r-cl(s_r)} + \mu_s(F \text{ is bad}) + \bar{\mu}(\bar{F} \text{ is bad}),$$

where $c = (f - \bar{h} - 2\epsilon)/\log 2 > 0$ and $F, \bar{F}$ denote the $s_r$-fillers in $X, \bar{X}$, respectively. It is easy to see that $EC^p$ is finite if $\sum Ec_{r+1}{}^p P(C \geq c_{r+1}) < \infty$, so it suffices to show that the three following series converge

$$\sum_{r=1}^{\infty} E(c_{r+1}{}^p 2^{m_r+n_r-cl(s_r)}), \quad \sum_{r=1}^{\infty} E(c_{r+1}{}^p \mu_s(F \text{ is bad})), \quad \sum_{r=1}^{\infty} E(c_{r+1}{}^p \bar{\mu}(\bar{F} \text{ is bad})).$$

Convergence is obtained by a careful study of analytic properties of the generating function of the filler length, the use of the Bernstein inequality to give exponential bounds for appropriate large deviation events, and repeated use of Hölder inequality.

The second step is a Bernoulli-to-Markov coding. Now $X$ is Bernoulli and $\bar{X}$ is mixing Markov with $h(\bar{X}) = \bar{h} < h = h(X)$. Moreover, as in [1], by extending $\bar{X}$ to another mixing Markov process with a slightly larger entropy (the extension is a coding of length one) we may assume that there exist a marker $\bar{M}$ in $\bar{X}$ and a single-word marker $M$ in $X$ such that the corresponding marker processes have the same distribution. Therefore the two marker processes can be identified as a common factor $\hat{X}$ of $X$ and $\bar{X}$. The bad fillers in $X$ and $\bar{X}$ are defined as in the Markov-to-Bernoulli case and the finiteness of $EC^p$ is concluded similarly by studying the three series (with some additional difficulties caused be the fact that the image is a Markov and not independent process).

**Remark 1.** The above methods allowed to compute moments of order $p$ with $p < 2$ for Markov-to-Bernoulli and Bernoulli-to-Markov coding, and left an open question whether the variance of the coding is finite. In a yet unpublished manuscript [7], the authors claim that there exists a universal finitary coding between Bernoulli schemes of unequal entropies which has exponential tails, i.e. $\text{Prob}(C > n)$ decays exponentially as $n \to \infty$, a condition which clearly implies that moments of all orders are finite. A code is universal in the following sense: if $A$ and $B$ are two alphabets, and a Bernoulli scheme on $B$ is given, of entropy $\bar{h}$, and in addition an $\epsilon > 0$ is given, then there exists a measurable subset of $A^{\mathbb{Z}}$ and a mapping $\phi$ from that set to $B^{\mathbb{Z}}$ such that if any Bernoulli scheme on $A$ is given, of entropy larger than $\bar{h} + \epsilon$, then $\phi$ is finitary with exponential tails.

Let us note, however, that this result does not imply that the Keane-Smorodinsky code has finite variance.

## 2.2. Case of equal entropies

In [11] Keane and Smorodinsky improved upon the result of Ornstein and showed:

**Theorem 5.** *If $h(p) = h(q)$, then $BS(p)$ and $BS(q)$ are finitarily isomorphic.*

A question was immediately posed as to under what additional assumptions the expected code length could be finite. It was known already that in the case of Meshalkin's code, the average code length was infinite, as for that particular code the probability $\text{Prob}(C > n)$ is equal to the probability that a simple random walk remains positive after $n$ steps.

First general statement in this direction was made by Parry ([21]), who provided a class of isomorphisms which had infinite expected code lengths, and therefore showed that entropy alone was not an invariant for finitary isomorphism with fect. An obstruction which was discovered by Parry involved the so-called information cocycle. Let us suppose that $\alpha = \{A_1, \ldots, A_k\}$ is a time-zero partition (also called state partition) of the shift space $(X, S, \mu)$, that is, $A_i = \{x : x_0 = i\}$. Let $\alpha^-$ denote the smallest $\sigma$-algebra containing $\bigvee_{i=1}^{\infty} S^{-i}\alpha$, and define the *information cocycle* of $S$ to be

$$I_S = I(\alpha|\ \alpha^-) = -\sum_i \chi_{A_i} \log \mu(A_i|\ \alpha^-).$$

In [21] Parry proved:

**Theorem 6.** *If $S$ and $T$ are finite state processes and if $\phi$ is a finitary isomorphism between $S$ and $T$ such that $\phi$ and $\phi^{-1}$ have finite expected code lengths, then the information cocycles $I_S$ and $I_T \circ \phi$ are cohomologous, with a finite valued and measurable transfer function, i.e. $I_S = I_T \circ \phi + g \circ S - g$.*

Note that here the dynamical systems were not assumed to be Markov.

It was then shown by Parry that particular dynamical systems which were known to be finitarily isomorphic by the results of Keane and Smorodinsky, had non-cohomologous information cocycles and therefore could not be finitarily isomorphic with fect. Among the above was Meshalkin's example of Bernoulli schemes $BS(\frac{1}{2}, \frac{1}{8}, \frac{1}{8}, \frac{1}{8}, \frac{1}{8})$ and $BS(\frac{1}{4}, \frac{1}{4}, \frac{1}{4}, \frac{1}{4})$ and some equal entropies Markov processes. First attempt in forming a set of invariants for isomorphism with fect between Markov shifts, was undertaken by Krieger. He (in [15]) defined, for a Markov shift with transition matrix $P$, a multiplicative subgroup $\Delta_P$ of positive reals in the following way:

$$\Delta_P = \{\frac{P(i_0, i_1) \cdots P(i_{n-1}, i_0)}{P(i_0, j_1) \cdots P(j_{n-1}, i_0)}\},$$

which are ratios of weights of cycles of equal lengths starting and ending in the same state. Krieger was then able to prove that if Markov shifts $P$ and $Q$ were finitarily isomorphic with fect then their respective delta groups $\Delta_P$ and $\Delta_Q$ were equal. He also gave examples of shifts with equal delta groups which could not be finitarily isomorphic with fect.

In 1981 Tuncel ([33]) introduced the so-called $\beta$-function, as the spectral radius of the matrix $P^t$ (where $P^t(i, j) = (P(i, j))^t$), and showed that this was an invariant of *regular* isomorphism. We shall not go into a detailed study of the classification up to a regular isomorphism; let us only note that this is a weaker notion than finitary isomorphism; an isomorphism $\phi$ is regular if both $\phi$ and its inverse $\phi^{-1}$ have bounded anticipation but can have infinite memory. Later Schmidt ([28], see also [29])improved on some of the results of Tuncel, showing in particular that $\beta$-function was an invariant for finitary isomorphism with fect.

In 1984 Parry and Schmidt ([22]) extended the notion of $\Delta_P$-group to that the $\Gamma_P$-group, generated by all weights $P(i_0, i_1) \cdots P(i_{n-1}, i_0)$ of cycles. They showed that for an aperiodic transition matrix $P$, the quotient group $\Gamma_P/\Delta_P$ is cyclic with a *distinguished* generator $c_P \Delta_P$. The main statement of [22] is that $\Gamma_P, \Delta_P$ and $c_P \Delta_P$ are invariants for finitary isomorphism with fect, which together with the previous result of Schmidt allowed to establish the following conjecture.

**Conjecture 1.** The quadruple $(\Gamma_P, \Delta_P, c_P \Delta_P, \beta_P)$ forms a complete set of invariants for finitary isomorphism with fect.

The above has been open for twenty years now, and the most general statement seems to be a recent result of Mouat and Tuncel:

**Theorem 7 ([19]).** *Let $P$ and $Q$ be primitive, stochastic matrices of the underlying Markov shifts, with the same $\Gamma, \Delta, c\Delta$ and $\beta$ invariants. If there exist states $I_0, J_0$ of the $P$-shift and the $Q$-shift, respectively, and there exist a nontrivial column vector $\mathbf{v}_r$ and a nontrivial row vector $\mathbf{v}_l$ such that $(P^n \mathbf{v}_r)(I_0) = (\mathbf{v}_l Q^n)(J_0)$ for all $n \geq 1$, then the two Markov shifts are finitarily isomorphic with fect.*

## 3. Applications

In this section we discuss a number of results which are closely related to finitary coding.

### 3.1. *m-dependent processes*

The marker method of [11] and [12] was extended by Smorodinsky ([31]) to prove that $m$-dependent processes of equal entropy were finitarily isomorphic. Let us recall here that a stationary process is called *m-dependent* if its past and future become independent, if separated by $m$ units of time. It is an easy observation that processes which are finite factors of independent processes are $m$-dependent, so Smorodinsky's result implies that equal entropy, finite factors of independent processes are finitarily isomorphic. It is natural to ask about finitary instead of finite factors, and this is stated in [31] as a conjecture:

**Conjecture 2 (Finitary factors conjecture).** Equal entropy, finitary factors of independent processes are finitarily isomorphic.

Let us keep in mind a well-known fact from Ornstein's theory that measurable factors of independent processes are (measurably) isomorphic, if they have the same entropy.

### 3.2. $\mathbb{Z}^d$-actions

Another direction in which the results of [10] have been extended is the action of $\mathbb{Z}^d$, $d \geq 2$, rather than the action of a single shift on $\mathbb{Z}$. In [9] del Junco considered two random fields on $\mathbb{Z}^2$: an ergodic Markov field $X$ and an independent process $Y$ such that $h(X) > h(Y)$. He then proved the existence of a finitary homomorphism from $X$ to $Y$. It was left as an open question whether a mixing Markov field was always a finitary factor of a Bernoulli process of higher entropy. Later ([2]) an example was given by van den Berg and Steif, of a Markov field which was not a finitary factor of any independent process, so the finitary factors conjecture fails in $d$ dimensions, $d \geq 2$.

Among the difficulties which arise when one considers $\mathbb{Z}^2$ actions instead of a $\mathbb{Z}$ action is the choice of markers. Let us recall that markers were nonoverlapping blocks; a marker in $\mathbb{Z}^2$ should therefore be a block which does not overlap itself under translation in any given direction, and that is a hard condition to fulfill. del Junco considers a multidimensional version of the Rokhlin tower lemma in order to define configurations in $\mathbb{Z}^2$ which have a number of disjoint shifts and which, together with a specific number of its shifts, almost cover the whole space; those configurations depend upon the multiple occurrences of blocks called markers. Skeletons of all ranks are defined in a highly nontrivial way, subsequently del Junco adapts and modifies when necessary the main ideas of [10], including the marriage lemma, to build the desired finitary coding between the random fields.

We complete this section by mentioning that it is the subject of the current research of the author of this note to show that the average *code volume* for the finitary coding from a Markov random field into a Bernoulli field of strictly smaller entropy, is finite. The generalization of the notion of code length to that of code volume in $\mathbb{Z}^2$ is straightforward. We do, however, propose an alternative approach to that of del Junco. A marker is a fixed block of a low probability of occurrence; when a marker occurs at a fixed coordinate, we consider a skeleton at this coordinate as a set of coordinates for which this marker occurrence is the closest, in the $L_1$-distance. This procedure gives rise to the partition of $\mathbb{Z}^2$, into the so-called Voronoi regions; the details will appear elsewhere.

A result of different flavor was obtained by Steif in [32], where he considered the so-called $T, T^{-1}$ process (also known as Random Walk in Random Scenery) in $d$-dimensional integer lattice (for definitions see [32] and references therein). A classical result of Kalikow is that the $T, T^{-1}$ process on $\mathbb{Z}$ associated to a simple random walk is not Bernoulli; Steif proves that in $\mathbb{Z}^d$ the second coordinate of this process is not a finitary factor of an independent process, he also applies this to study some properties of the Ising model in statistical mechanics.

### 3.3. Countable state processes

It turns out that the assumption that a Markov shift be a <u>finite</u> state process, is an essential one, as far as classification up to a finitary isomorphism is concerned. Smorodinsky proved (unpublished, see e.g. [16]) that an ergodic automorphism of compact abelian group can only be finitarily Bernoulli (that means: finitarily isomorphic to a Bernoulli shift) if it is *exponentially recurrent* (i.e. if $U$ is an open set and $r_U$ is the first-return time function, then $\text{Prob}(r_U > n)$ decays exponentially fast with $n$). Smorodinsky used this result to construct a countable state Markov shift which is measurably isomorphic to a Bernoulli shift but has polynomially decaying return times, hence cannot be finitarily isomorphic to a Bernoulli scheme. Lind ([16]) went on in that direction to prove that ergodic automorphisms of compact abelian groups are exponentially recurrent, and that was a step forward in trying to resolve a question as to whether ergodic group automorphisms are finitarily isomorphic to Bernoulli shifts. Let us recall here that some special classes, like hyperbolic toral automorphisms, are known to be finitarily Bernoulli.

Rudolph ([27]) completed Smorodinsky's work in showing that a countable state, mixing Markov shift of finite entropy is finitarily Bernoulli if and only if the chain has exponentially decaying return times; in particular, it follows that the so-called $\beta$-automorphisms are finitarily isomorphic to independent processes of the same finite entropy. We wish to note that this result depends heavily upon the characterization of the processes finitarily isomorphic to Bernoulli shifts ([26]), a construction which is a very involved generalization of [11]. More recently, Keane and Steif ([13]) proved that $T, T^{-1}$ process associated to a 1-dimensional random walk with positive drift is finitarily isomorphic to an independent process, using an intermediate countable state Markov shift, and results of [27].

Finally, we remark that Petit ([25]) extended the methods of Keane and Smorodinsky ([11]) to show that two infinite entropy Bernoulli schemes on countable state space are finitarily isomorphic.

Let us close this paragraph by mentioning that despite all developments, the classification of countable state Markov processes up to a finitary isomorphism remains an intricate task, and there seems to be a need for methods which would both be more elementary and more easily understood than the present ones.